\documentclass[a4paper]{article}

\usepackage[english]{babel}
\usepackage[utf8x]{inputenc}
\usepackage[T1]{fontenc}

\usepackage[top=3cm,bottom=2cm,left=3cm,right=3cm,marginparwidth=1.75cm]{geometry}
\usepackage{authblk}
\usepackage{blindtext}
\usepackage{amsmath, amsthm, amssymb, esint, mathtools, mathabx}
\usepackage{amsfonts}
\usepackage{graphicx}
\graphicspath{{/texfiles/}}
\usepackage[colorinlistoftodos]{todonotes}
\usepackage[colorlinks=true, allcolors=blue]{hyperref}
\usepackage{multicol}

\numberwithin{equation}{section}

\title{Anti-uniformity norms, anti-uniformity functions and their algebras on Euclidean spaces}
\author{A. Martina Neuman}
\affil{Department of Mathematics, New York University, Shanghai}
\affil[ ]{\textit{marsneuman@nyu.edu}}

\begin{document}

\maketitle

\begin{abstract} 
Let $k\geq 2$ be an integer. Given a uniform function $f$ - one that satisfies $\|f\|_{U(k)}<\infty$, there is an associated anti-uniform function $g$ - one that satisfied $\|g\|_{U(k)}^{*}$. The question is, can one approximate $g$ with the Gowers-Host-Kra dual function $D_{k}f$ of $f$? Moreover, given the generalized cubic convolution products $D_{k}(f_{\alpha}:\alpha\in\tilde{V}_{k})$, what sorts of algebras can they form?\\
In short, this paper explores possible structures of anti-uniformity on Euclidean spaces.
\end{abstract}


\tableofcontents

\section{Introduction}

\noindent
This exposition's main aim is not to prove one central result, but to produce a series of discoveries around a main theme: to understand the anti-uniform functions and the related cubic convolution products as well as algebraic structures that can be formed from them. The anti-uniform structures here are taken from the Gowers-Host-Kra structures ([3], [5], [7], [6], [12]). As the Gowers-Host-Kra structure is a complicated structure that needs some build-ups before making any coherent statement, an impatient reader can skip to the next sections first before going back to this section for definition. The main results are in {\bf Section 4} and {\bf Section 5}. See also the {\bf Acknowledgement} below for more initial motivation of this project.\\

\noindent
In this exposition, $k\geq 2$ and all the functions involved are nonnegative, except for in this initial setup, where structures are defined for general real-valued functions and for in {\bf Section 5}, where algebras of functions are considered.\\

\noindent
Fix a dimension $d$. Then given a measurable function $f:\mathbb{R}^{d}\to\mathbb{R}$ and an integer $k\geq 1$, one defines the $k$th Gowers-Host-Kra norm of $f$, $\|f\|_{U(k)}$ inductively as follows,
\begin{equation}\label{eq 1.1}\|f\|_{U(1)} = |\int_{\mathbb{R}^{d}} f\,dx| \,\,\,\text{ and }\,\,\,\|f\|_{U(k+1)}^{2^{k+1}} = \int_{\mathbb{R}^{d}} \|f^{h}\cdot f\|_{U(k)}^{2^{k}}\,dh \,\,\text{, if } k\geq 1.\end{equation}
Here $f^{h}(x) = f(x+h)$. It's well-known that for $k\geq 2$, $\|\cdot\|_{U(k)}$ is a proper norm and that $\|f\|_{U(2)} = \|\hat{f}\|_4$; see [5], [12]. It's also known that for every $k\geq 1$, there exist $A(k,d) = A(k)$ and $p_{k}=2^{k}/(k+1)$ such that 
\begin{equation} \label{eq 1.2} \|f\|_{U(k)}\leq A(k)\|f\|_{p_{k}}.\end{equation}

\noindent
A noninductive definition of the $k$th Gowers-Host-Kra norm can be given as follows. Let $V_{k} = \{0,1\}^{k}$ ($V$ stands for "vertices"). Let $\vec{h}=(h_1,\cdots,h_{k})\in (\mathbb{R}^{d})^{\otimes k}$. If $\alpha = (\alpha_1,\cdots,\alpha_{k})\in V_{k}$ then 
\begin{equation*}\alpha\cdot\vec{h}:= \alpha_1 h_1 +\cdots+\alpha_{k} h_{k}.\end{equation*}
Then one can confirm that for $k\geq 2$ (see [5], [12]), $\ref{eq 1.1}$ is the same as,
\begin{equation}\label{eq 1.3} \|f\|_{U(k)}^{2^{k}} = \int_{\mathbb{R}^{(k+1)d}} \prod_{\alpha\in V_{k}} f(x+\alpha\cdot\vec{h})\,dx d\vec{h}.\end{equation}

\noindent
The integral apparatus on the RHS of $\ref{eq 1.3}$ can be generalized for measurable $f_{\alpha}$, $\alpha\in V_{k}$, that are not necessarily the same, as follows,
\begin{equation}\label{eq 1.4} \mathcal{T}_{k}(f_{\alpha}:\alpha\in V_{k}):=\int_{\mathbb{R}^{(k+1)d}} \prod_{\alpha\in V_{k}} f_{\alpha}(x+\alpha\cdot\vec{h})\,dxd\vec{h}.\end{equation}
Temporarily name this the $k$th Gowers inner product of $(f_{\alpha}:\alpha\in V_{k})$.\\

\noindent
Some properties of $p_{k}, A(k)$ that are worth noticing, are:
\begin{itemize}
    \item $p_{k}\geq 1$ for all $k\geq 1$ and $p_{k}\uparrow\infty$ as $k\uparrow\infty$.
    \item $A(k)\leq 1$ for all $k\geq 1$ (and for all dimensions $d$); $A(k)$'s are increasing starting with $k\geq 3$.
\end{itemize}
There are more interesting properties of $p_{k}$ and $A(k)$ relating to the exponents and the optimal constants involved in the sharp Young's inequality. See [5], [12]. One more useful knowledge is, [12],\\

\noindent
{\bf Cauchy-Schwarz-Gowers inequality:} 
\begin{equation}\label{eq 1.5} \mathcal{T}_{k}(f_{\alpha}:\alpha\in V_{k}) \leq\prod_{\alpha\in V_{k}}\|f_{\alpha}\|_{U(k)}\leq A(k)^{2^{k}}\prod_{\alpha\in V_{k}}\|f_{\alpha}\|_{p_{k}}.\end{equation}
If $f_{\alpha}=f$ for all $\alpha\in V_{k}$ then the first inequality of $\ref{eq 1.5}$ becomes $\ref{eq 1.1}$ and the second becomes $\ref{eq 1.2}$.\\

\noindent
{\bf Dual functions as cubic convolution product:} Let $\tilde{V}_{k}=V_{k}\setminus\{\vec{0}\}$. Define the $k$th Gowers-Host-Kra dual function of $f$ as,
\begin{equation*} D_{k}f(x) = \int_{\mathbb{R}^{kd}}\prod_{\alpha\in\tilde{V}_{k}}f(x+\alpha\cdot\vec{h})\,d\vec{h}.\end{equation*} 
It's clear that $D_{k}f(x)$ is continuous if $f$ is locally integrable. There is a concept of generalized {\it cubic convolution product} of $(f_{\alpha}:\alpha\in\tilde{V}_{k})$,
\begin{equation*}D_{k}(f_{\alpha}:\alpha\in V_{k})(x) =\int_{\mathbb{R}^{kd}}\prod_{\alpha\in\tilde{V}_{k}}f_{\alpha}(x+\alpha\cdot\vec{h})\,d\vec{h}.\end{equation*}
Note that for every $t\in\mathbb{R}$, \begin{equation}\label{eq 1.6} D_{k}(tf)(x) = t^{2^{k}-1}D_{k}f(x).\end{equation}\\

\noindent
{\bf The exponents:} There are main Lebesgue integration exponents used here, $p_{k}, q_{k}, s_{k}$. The first has been mentioned above and the latter two will be later when the need comes. It's worth to give their definitions again in a central place, however. Let $k\geq 2$ be an integer, then
\begin{align*}
p_{k} &=\frac{2^{k}}{k+1}\\
q_{k} &=\frac{2^{k}-1}{k}\\
s_{k} &=\frac{2^{k}}{2^{k}-k-1}.\end{align*}
In other words, $s_{k}$ is the H\"{o}lder's conjugate of $p_{k}$.

\subsection{Acknowledgement}

\noindent
The author would like to thank Professor Bryna Kra for her suggestion of this project. It starts out as an effort to find an analogue of {\bf Theorem 3.8} and {\bf Corollary 3.11} in [8] in the Euclidean settings. The author decided to find more possible Euclidean parallels she could find with the other ergodic results in the said paper.\\

\noindent
Apart from the parallels, there's quite a contrast between the two worlds. For instance, due to the noncompact nature of $\mathbb{R}^{d}$ as well as the strong connection between the sharp Young's convolution inequality and the Gowers inequalities (the sharp constants as well as the Lebesgue exponents), {\bf Lemma 2.3} and {\bf Lemma 3.4} don't occur exactly in the Euclidean settings. The analogues are  {\bf Lemma 1} and {\bf Lemma 2}, respectively, where the exponents involved are not the same as those in the ergodic setting, and they are not even equal, as in the ergodic setting. As a result, there is no connection between the algebra norms $\|\cdot\|_{\mathbb{A}(k)},\|\cdot\|_{\mathbb{F}(k)}$ with the anti-uniform norm $\|\cdot\|_{U(k)}^{*}$ in Euclidean spaces, unlike in the ergodic setting. See {\bf Section 2}, {\bf Section 3}, {\bf Section 5} for the definitions of these norms.\\

\noindent
Due to the noncompact nature of the Euclidean spaces, the shrinking of the uniform spaces $\mathcal{U}(k)$ and the corresponding expansion of the anti-uniform spaces $\mathcal{U}^{*}(k)$ as $k\uparrow\infty$ are not observed. Hence, analogues for {\bf Theorem 3.8} and {\bf Corollary 3.11} in [8], which are {\bf Theorem 4} and {\bf Corollary 5}, respectively, can go as far as giving an approximation by an anti-uniform function a exactly level $k$ and not at any other levels. Hence, there's no analogue of {\bf Theorem 3.9} in [8].\\

\noindent
Lastly, the terminologies and the notations used here are vastly borrowed from [8]. Some divergences in the notations are due to early clashes with other notations that the author has used long before. 

\subsection{Caution}

\noindent
Since an integral that defines a Gowers-Host-Kra structure typically involves integration over $\mathbb{R}^{(k+1)d}$ or $\mathbb{R}^{kd}$, writing down such integral in a neat manner is quite a challenge. The author adopts the convention that, often, only one integral sign is written down, and the indication of the integration variables will be given. For example, $\int \cdots \,d\vec{t}d\vec{h}$, where $\vec{t},\vec{h}\in\mathbb{R}^{M}$, indicates an integration over $\mathbb{R}^{M}\times\mathbb{R}^{M}$.\\

\noindent
With all the constants having been introduced, in what follows, the implicit calculated constants are only dependent on the index $k$ and the dimension $d$.

\section{Cubic convolution product}

\noindent 
If $\beta\in V_{k}$ then one understands the notation $(\beta,\delta)$ as indicating a member of $V_{k+1}$ that begins in $\beta$ and ends in $\delta\in\{0,1\}$.\\

\noindent
Let $q_{k} = \frac{2^{k}-1}{k}$. One has the following domination,\\

\noindent
{\bf Lemma 1.} For every $k\geq 2$, 
\begin{equation} \label{eq 2.1} D_{k}(f_{\alpha}:\alpha\in \tilde{V}_{k})(x)\lesssim\prod_{\alpha\in\tilde{V}_{k}}\|f_{\alpha}\|_{q_{k}}.\end{equation}
{\it Proof:} The proof is based on induction on $k$. Let $k=2$. Then $\alpha=\{01,10,11\}$ if $\alpha\in\tilde{V}_2$, and, \begin{equation*}D_2(f_{01},f_{10},f_{11})(x)=\int f_{01}(x+t_1)f_{10}(x+t_2)f_{11}(x+t_1+t_2)\,dt_1\,dt_2.\end{equation*}
Denote $g(x+\cdot) = g^{x}(\cdot)$ and $\tilde{g}(x)=g(-x)$. Then
\begin{multline*}D_2(f_{01},f_{10},f_{11})(x) = \int f^{x}_{01}(t_1)f^{x}_{10}(t_2)f^{x}_{11}(t_1+t_2)\,dt_1 dt_2\\ = \int f^{x}_{01}\cdot(\tilde{f}^{x}_{10}\ast f^{x}_{11})(t_1)\,dt_1\lesssim\|f_{11}\|_{3/2}\|\tilde{f}^{x}_{10}\ast f^{x}_{11}\|_{3}\lesssim \|f_{11}\|_{3/2}\|f_{10}\|_{3/2}\|f_{01}\|_{3/2}\end{multline*}
where the first and last inequalities follow from H\"{o}lder's and Young's convolution inequalities, respectively. \\

\noindent
Now assume $\ref{eq 2.1}$ is true for some $k\geq 2$. If $\alpha\in\tilde{V}_{k+1}$ then $\alpha = (\beta,1), (\beta,0), (\vec{0},1)$ for some $\beta\in\tilde{V}_{k}$. Then,
\begin{equation*}D_{k+1}(f_{\alpha}:\alpha\in\tilde{V}_{k+1})(x) = \int \prod_{\beta\in\tilde{V}_{k}}f_{(\beta,0)}(x+\beta\cdot\vec{t})f_{(\beta,1)}(x+\beta\cdot\vec{t}+h)f_{(\vec{0},1)}(x+h)\,d\vec{t}dh.\end{equation*}
Integrating this last integral expression first in $\vec{t}\in\mathbb{R}^{k}$, we can use the induction hypothesis to obtain,
\begin{multline*}D_{k+1}(f_{\alpha}:\alpha\in\tilde{V}_{k+1})(x)\lesssim\int f_{(\vec{0},1)}(x+h)\prod_{\beta\in\tilde{V}_{k}}\|f^{h}_{(\beta,1)}\cdot f_{(\beta,0)}\|_{q_{k}}\,dh\\\lesssim\|f_{(\vec{0},1)}\|_{q_{k+1}}\prod_{\beta\in\tilde{V}_{k}}\left(\int \|f^{h}_{(\beta,1)}\cdot f_{(\beta,0)}\|_{q_{k}}^{r}\,dh\right)^{1/r} .\end{multline*}
The second inequality above follows from H\"{o}lder's inequality with the exponent $r$ satisfying $q_{k+1} + (2^{k}-1)\cdot\frac{1}{r}=1$. Now
\begin{equation*}\int \|f^{h}_{(\beta,1)}\cdot f_{(\beta,0)}\|_{q_{k}}^{r}\,dh = \int\left(\int f^{q_{k}}_{(\beta,1)}(x+h)f^{q_{k}}_{(\beta,0)}(x)\,dx \right)^{r/q_{k}}\,dh = \|f^{q_{k}}_{(\beta,1)}\ast \tilde{f}^{q_{k}}_{(\beta,0)}\|_{r/q_{k}}^{r/q_{k}}.\end{equation*}
Applying Young's convolution inequality one more time to the last norm above, one has, 
\begin{equation*}\int \|f^{h}_{(\beta,1)}\cdot f_{(\beta,0)}\|_{q_{k}}^{r}\,dh = \|f^{q_{k}}_{(\beta,1)}\ast \tilde{f}^{q_{k}}_{(\beta,0)}\|_{r/q_{k}}^{r/q_{k}}\lesssim \left(\|f^{q_{k}}_{(\beta,1)}\|_{s}\|f^{q_{k}}_{(\beta,0)}\|_{s}\right)^{r/q_{k}}\end{equation*}
with the exponent $s$ satisfying $1+q_{k}/r=2/s$. Simple arithmetic gives, $\frac{1}{q_{k}s} = \frac{k+1}{2^{k+1}-1}$, which is the desired exponent $q_{k+1}$. Putting all these back together, we obtain $\ref{eq 2.1}$ for the case of $k+1$,
\begin{equation*}D_{k+1}(f_{\alpha}:\alpha\in\tilde{V}_{k+1})(x)\lesssim\prod_{\alpha\in\tilde{V}_{k+1}}\|f_{\alpha}\|_{q_{k+1}}.\end{equation*} 
This concludes the induction step. $\qed$\\

\noindent
Recall that $\ref{eq 1.2}$ says, $\|f\|_{U(k)}\lesssim\|f\|_{p_{k}}$, which means one can consider the space $L^{p_{k}}$ to be endowed with the uniform norm $U(k)$. Denote this space, $\mathcal{U}(k)=(L^{p_{k}},\|\cdot\|_{U(k)})$. Then its dual space can be considered as subspace of $L^{s_{k}}$ with $s_{k} = (1-p_{k}^{-1})^{-1} = \frac{2^{k}}{2^{k}-k-1}$. Denote the dual space as $\mathcal{U}^{*}(k)$. Call this dual space the $k$th anti-uniform space, and the functions in this space, the $k$th anti-uniform functions. The $\|\cdot\|_{U(k)}^{*}$ is given by duality - call it the $k$th anti-uniform norm. More specifically, $g\in\mathcal{U}^{*}(k)$ if $g\in L^{s_{k}}$ and if
\begin{equation*}\sup\{|\langle g,f\rangle|: f\in L^{p_{k}}, \|f\|_{U(k)}\leq 1\}<\infty.\end{equation*}
Then $\|g\|_{U(k)}^{*}$ is defined to be this supremum value. Note that in the supremum above, it's enough to consider only those $f\in L^{p_{k}}$ such that $\|f\|_{U(k)}=1$. For suppose that there is an optimizer $f$ such that $\|g\|^{*}_{U(k)}=: C =\langle g,f\rangle$ and $\delta=\|f\|_{U(k)}<1$. Then $\|f/\delta\|_{U(k)}=1$ and $\langle g,f/\delta\rangle =(\delta)^{-1}\langle g,f\rangle > C$, which gives a contradiction. It's also easy to see that for $f\in\mathcal{U}(k)$ and $g\in\mathcal{U}^{*}(k)$
\begin{equation}\label{eq 2.2} \|f\|_{U(k)}\|g\|_{U(k)}^{*}\geq\langle g,f\rangle.\end{equation}

\noindent
Now if $g\in\mathcal{U}^{*}(k)$ then $g\in L^{s_{k}}$. Take $f\in L^{p_{k}}$ such that $\|f\|_{p_{k}}= 1$ and $\|g\|_{s_{k}}=\langle g,f\rangle$. Then \begin{equation*}\|g\|^{*}_{U(k)}\geq\langle g,f\rangle = \|g\|_{s_{k}}.\end{equation*}
However, $\ref{eq 2.2}$ gives a better inequality, since if $\|f\|_{p_{k}}=1$ then $\|f\|_{U(k)}\leq A(k)$ by $\ref{eq 1.2}$, hence, \begin{equation*}A(k)\|g\|_{U(k)}^{*}\geq\|f\|_{U(k)}\|g\|_{U(k)}^{*}\geq\|g\|_{s_{k}}\end{equation*}
which is
\begin{equation}\label{eq 2.3}\|g\|^{*}_{U(k)}\geq A(k)^{-1}\|g\|_{s_{k}}.\end{equation}

\section{Anti-uniform functions}

\noindent
{\bf Lemma 2.} Let $f_{\alpha}\in L^{p_{k}}$ for every $\alpha\in\tilde{V}_{k}$. Then 
\begin{equation}\label{eq 3.1}
\|D_{k}(f_{\alpha}:\alpha\in\tilde{V}_{k})\|_{U(k)}^{*}\leq\prod_{\alpha\in\tilde{V}_{k}}\|f_{\alpha}\|_{U(k)}\lesssim\prod_{\alpha\in\tilde{V}_{k}}\|f_{\alpha}\|_{p_{k}}.\end{equation}
{\it Proof:} Let $g=D_{k}(f_{\alpha}:\alpha\in\tilde{V}_{k})$. Take $h\in\mathcal{U}(k)$ such that $\|h\|_{U(k)}\leq 1$. Then it follows from $\ref{eq 1.2}, \ref{eq 1.4}, \ref{eq 1.5}$ that
\begin{multline*}\langle g,h\rangle=\int h(x+\vec{0}\cdot\vec{t})\prod_{\alpha\in\tilde{V}_{k}}f_{\alpha}(x+\alpha\cdot\vec{t})\,d\vec{t}dx = \mathcal{T}_{k}(f_{\vec{0}}=h,f_{\alpha}:\alpha\in\tilde{V}_{k})\\ \leq\|h\|_{U(k)}\prod_{\alpha\in\tilde{V}_{k}}\|f_{\alpha}\|_{U(k)}\lesssim \|h\|_{U(k)}\prod_{\alpha\in\tilde{V}_{k}}\|f_{\alpha}\|_{p_{k}}\end{multline*}
which yields $\ref{eq 3.1}$. $\qed$\\

\noindent
On the other hand, from $\ref{eq 1.3}$ and $\ref{eq 2.2}$
\begin{equation*}\|f\|_{U(k)}^{2^{k}}=\langle D_{k}f,f\rangle\leq\|D_{k}f\|^{*}_{U(k)}\|f\|_{U(k)}\end{equation*}
which means $\|D_{k}f\|_{U(k)}^{*}\geq\|f\|_{U(k)}^{2^{k}-1}$. This and $\ref{eq 3.1}$ yield \begin{equation}\label{eq 3.2} \|D_{k}f\|_{U(k)}^{*}=\|f\|_{U(k)}^{2^{k}-1}.\end{equation}

\noindent
{\bf Lemma 3.} The unit ball of the $k$th anti-uniform space $\mathcal{U}^{*}(k)$ is the closed convex hull in $L^{s_{k}}$ of the set
\begin{equation*}K = \{D_{k}f: f\in L^{p_{k}}, \|f\|_{U(k)}\leq 1\}\end{equation*}
{\it Proof:} Let $B$ be the said unit ball. Denote $\bar{K}:=\overline{K}^{L^{s_{k}}}$. By $\ref{eq 3.2}$, if $\|f\|_{U(k)}\leq 1$ then $\|D_{k}f\|^{*}_{U(k)}\leq 1$ hence $D_{k}f\in B$. That means $K\subset B$. Now from $\ref{eq 2.3}$, $B$ is contained in the unit ball of the reflexive Banach space $L^{s_{k}}$ and is a weak$-*$ compact subset of this space. Hence $\overline{B}^{L^{s_{k}}}=B$, and so $\bar{K}\subset B$.\\

\noindent
Now suppose $\bar{K}\not\supset B$, then there exists $h\in B$ but $h\not\in\bar{K}$. That means $\|h\|^{*}_{U(k)}\leq 1$ but $h\not\in\bar{K}$. By the Hahn-Banach separation theorem (see [1]), that means there exists $f\in L^{p_{k}}$ such that
\begin{enumerate}
    \item $\langle f,g\rangle\leq 1$ for every $g\in\bar{K}$, yet,
    \item $\langle f,h\rangle > 1.$
\end{enumerate}
Through $\ref{eq 2.2}$ and the hypothesis, the second property implies $\|f\|_{U(k)} > 1$. Now take $F = \|f\|_{U(k)}^{-1}f$, then $\|F\|_{U(k)}=1$ and so, $D_{k}F\in K\subset\bar{K}$. From the first property and $\ref{eq 1.3}, \ref{eq 1.6}$, $\|f\|_{U(k)}= \|f\|_{U(k)}^{-2^{k}+1}\langle D_{k}f,f\rangle=\langle D_{k}F,f\rangle\leq 1$, which is a contradiction to the mentioned consequence of the second property. $\qed$

\section{Approximation results for anti-uniform functions}

\noindent
{\bf Theorem 4.} Let $k\geq 2$ be an integer. For every anti-uniform function $g$ with $\|g\|_{U(k)}^{*}=1$ and $\delta>0$, then $g$ can be written as 
\begin{equation} \label{eq 4.1} g=D_{k}F+H\end{equation} such that 
\begin{equation} \label{eq 4.2} \|F\|_{p_{k}}\leq 1/\delta\,\,\, \text{ and }\,\,\,\|F\|_{U(k)}\leq 1\,\,\,\text{ and }\,\,\,\|H\|_{s_{k}}\leq\delta.\end{equation}
{\it Proof:} Fix such a $k$ and $\delta>0$. For $f\in L^{p_{k}}$, define \begin{equation}\label{eq 4.3}
\vvvert f\vvvert =(\|f\|_{U(k)}^{2^{k}}+\delta^{2^{k+1}}\|f\|_{p_{k}}^{2^{k}})^{1/2^{k}}.\end{equation}
From $\ref{eq 1.2}$, $\vvvert f\vvvert\lesssim_{\delta}\|f\|_{p_{k}}$; that means this norm is well-defined on $L^{p_{k}}$. On the other hand, clearly $\vvvert f\vvvert\gtrsim_{\delta}\|f\|_{p_{k}}$. Hence the norm $\vvvert\cdot\vvvert$ is equivalent to the norm $\|\cdot\|_{p_{k}}$ on $L^{p_{k}}$. Consequently, $(L^{p_{k}},\vvvert\cdot\vvvert)$ is a reflexive Banach space. Now define a dual norm $\vvvert\cdot\vvvert^{*}$ on $L^{s_{k}}$ by,
\begin{equation}\label{eq 4.4} \vvvert g\vvvert^{*}=\sup\{|\langle f,g\rangle|: f\in L^{p_{k}}, \vvvert f\vvvert\leq 1\}.\end{equation}
Then from duality, this dual norm is equivalent to $\|\cdot\|_{s_{k}}$. That means, $(L^{p_{k}},\vvvert\cdot\vvvert)^{*}=(L^{s_{k}},\vvvert\cdot\vvvert^{*})$. Since $\vvvert f\vvvert\geq\|f\|_{U(k)}$, by a simple argument as with the one used in $\ref{eq 2.3}$, one has that $\vvvert g\vvvert^{*}\leq\|g\|^{*}_{U(k)}$. \\

\noindent
Fix $g$ with $\|g\|^{*}_{U(k)}=1$ and let $C=\vvvert g\vvvert^{*}\leq\|g\|^{*}_{U(k)}\leq 1$. Let $\mathcal{G}=C^{-1}g$ so that $\vvvert\mathcal{G}\vvvert^{*}=1$.\\
From $\ref{eq 4.4}$ and reflexivity, there exists $\mathcal{F}\in L^{p_{k}}$ such that \begin{equation}\label{eq 4.5}\vvvert\mathcal{F}\vvvert=1\,\,\,\text{ and }\,\,\,\langle\mathcal{G},\mathcal{F}\rangle = 1.\end{equation}
This means, from $\ref{eq 4.3}$, $\|\mathcal{F}\|_{U(k)}\leq 1$ and $\|\mathcal{F}\|_{p_{k}}\leq 1/\delta$. Now for every $\phi\in L^{p_{k}}$ and $t\in\mathbb{R}$,
\begin{equation}\label{eq 4.6}\|\mathcal{F}+t\phi\|_{U(k)}^{2^{k}}=\langle\mathcal{F}+t\phi, D_{k}(\mathcal{F}+t\phi)\rangle = \|\mathcal{F}\|_{U(k)}^{2^{k}}+2^{k}t\langle D_{k}\mathcal{F},\phi\rangle + o(t).\end{equation}
The last expression of $\ref{eq 4.6}$ is obtained through simple expansion, changing of integration variables and collecting like-terms. The notation $o(t)$ means a collection of terms whose absolute value is at most $o(t)$, which satisfies $\lim_{t\to 0} o(t)/t=0$. On the other hand, one has
\begin{equation}\label{eq 4.7}\|\mathcal{F}+t\phi\|_{p_{k}}^{2^{k}}=\|\mathcal{F}\|_{p_{k}}^{2^{k}}+2^{k}t\|\mathcal{F}\|_{p_{k}}^{2^{k}-1}\langle\mathcal{F}^{p_{k}-1},\phi\rangle + o(t).\end{equation}
Now $\ref{eq 4.7}$ follows from that $p_{k}=\frac{2^{k}}{k+1}$ and the Taylor expansion of $T(t) = \|\mathcal{F}+t\phi\|_{p_{k}}^{2^{k}}$ at $t=0$. Then $\ref{eq 4.3}, \ref{eq 4.6},\ref{eq 4.7}$ altogether give,
\begin{multline*} \vvvert\mathcal{F}+t\phi\vvvert^{2^{k}}=\|\mathcal{F}+t\phi\|_{U(k)}^{2^{k}}+\delta^{2^{k+1}}\|\mathcal{F}+t\phi\|_{p_{k}}^{2^{k}} \\ = \vvvert\mathcal{F}\vvvert^{2^{k}}+ 2^{k}t\langle D_{k}\mathcal{F},\phi\rangle+ \delta^{2^{k+1}}2^{k}t\|\mathcal{F}\|_{p_{k}}^{2^{k}-1}\langle \mathcal{F}^{p_{k}-1},\phi\rangle + o(t),\end{multline*}
which when raised to the power $1/2^{k}$ - and keeping in mind $\ref{eq 4.5}$ - gives
\begin{equation}\label{eq 4.8} \vvvert\mathcal{F}+t\phi\vvvert = 1 + t\langle D_{l}f',\phi\rangle +\delta^{2^{k+1}}t\|\mathcal{F}\|_{p_{k}}^{2^{k}-1}\langle \mathcal{F}^{p_{k}-1},\phi\rangle+o(t).\end{equation}
$\ref{eq 4.8}$ is another application of Taylor expansion. Now from $\ref{eq 4.5}$ and duality, for every $\phi\in L^{p_{k}}$ one has, 
\begin{equation*}1 + t\langle\mathcal{G},\phi\rangle = \langle\mathcal{F},\mathcal{G}\rangle + t\langle\mathcal{G},\phi\rangle= \langle\mathcal{G},\mathcal{F}+t\phi\rangle\leq \vvvert\mathcal{F}+t\phi\vvvert\end{equation*}
which together with $\ref{eq 4.8}$ gives, 
\begin{equation*}
1+t\langle\mathcal{G},\phi\rangle\leq 1 + t\langle D_{k}\mathcal{F},\phi\rangle +\delta^{2^{k+1}}t\|\mathcal{F}\|_{p_{k}}^{2^{k}-1}\langle \mathcal{F}^{p_{k}-1},\phi\rangle+o(t)\end{equation*}
which in turns implies,
\begin{equation*}
\langle\mathcal{G},\phi\rangle\leq \langle D_{k}\mathcal{F},\phi\rangle+\delta^{2^{k+1}}\|\mathcal{F}\|_{p_{k}}^{2^{k}-1}\langle\mathcal{F}^{p_{k}-1},\phi\rangle .\end{equation*}
Since this holds for every $\phi\in L^{p_{k}}$, by duality, one has, 
\begin{equation*}\mathcal{G}= D_{k}\mathcal{F}+\delta^{2^{k+1}}\|\mathcal{F}\|_{p_{k}}^{2^{k}-1}\mathcal{F}^{p_{k}-1}.\end{equation*}
Set $F = C^{1/(2^{k}-1)}\mathcal{F}$ and $H=C\delta^{2^{k+1}}\|\mathcal{F}\|_{p_{k}}^{2^{k}-1}\mathcal{F}^{p_{k}-1}$, then
\begin{equation*}g= C\mathcal{G} = D_{k}F+ H.\end{equation*}
Now $\|F\|_{U(k)}\leq 1$ and $\|F\|_{p_{k}}\leq 1/\delta$, since $\|\mathcal{F}\|_{U(k)}\leq 1$ and $C\leq 1$, and \begin{equation*}\|H\|_{s_{k}} = C\delta^{2^{k+1}}\|\mathcal{F}\|_{p_{k}}^{2^{k}-1}\|\mathcal{F}^{p_{k}-1}\|_{s_{k}} = C\delta^{2^{k+1}}\|\mathcal{F}\|_{p_{k}}^{2^{k}-1}\|\mathcal{F}\|_{p_{k}}^{p_{k}/s_{k}}\leq\delta\end{equation*}
which are all the requirements of $\ref{eq 4.1}, \ref{eq 4.2}$. $\qed$\\

\noindent
{\bf Corollary 5.} Let $k\geq 2$ be an integer. Let $\phi$ be such that $\|\phi\|_{p_{k}}\leq 1$ $\|\phi\|_{U(k)}=\theta >0$. Then there exists a function $f$ such that \begin{equation}\label{eq 4.9} \|f\|_{p_{k}}\leq 1\,\,\,\text{ and }\,\,\,\langle D_{k}f,\phi\rangle >(\theta/2)^{2^{k}}.\end{equation}
{\it Proof:} Given such $h$, there exists $g$ with $\|g\|_{U(k)}^{*}=1$ and $\langle g,\phi\rangle>\theta$. From {\bf Theorem 4}, for every $\delta>0$, one has, $g= D_{k}F+H$ with $\|F\|_{p_{k}}\leq 1/\delta$ and $\|H\|_{s_{k}}\leq\delta$. That means $\langle H,\phi\rangle\leq\|\phi\|_{p_{k}}\|H\|_{s_{k}}\leq\delta$ and $\langle D_{k}F,\phi\rangle >\theta-\delta$. Let $f=\delta F$ and take $\delta=\theta/2$, one arrives at $\ref{eq 4.9}$. $\qed$

\section{The associated algebras}

\noindent
In this section, functions are allowed to have negative values.\\

\noindent
A similar argument as in $\ref{eq 3.1}$ (or use $\ref{eq 1.5}$) gives that, if $g=D_{k}(f_{\alpha}:\alpha\in\tilde{V}_{k})$ with $f_{\alpha}\in L^{p_{k}}$, then,
\begin{equation*}\|g\|_{s_{k}}\lesssim\prod_{\alpha\in\tilde{V}_{k}}\|f_{\alpha}\|_{p_{k}}.\end{equation*}

\noindent
Follow from this and {\bf Lemma 1}, there are two algebras associated with the cubic convolution product functions.

\subsection{The algebra $\mathbb{A}(k)$}

\noindent
{\bf Definition 6.} Let $k\geq 2$ be an integer. Define $\mathcal{K}(k)$ to be the convex hull of 
\begin{equation*}\{D_{k}(f_{\alpha}:\alpha\in\tilde{V}_{k}): f_{\alpha}\in L^{q_{k}}\text{ and } \prod_{\alpha\in\tilde{V}_{k}}\|f_{\alpha}\|_{q_{k}}\leq 1\}\end{equation*}

\noindent
{\bf Lemma 7.} If $f_{\alpha}\in L^{q_{k}}$, $\alpha\in\tilde{V}_{k}$, then $D_{k}(f_{\alpha}:\alpha\in\tilde{V}_{k})(x)$ is a continuous function.\\
{\it Proof:} If $f_{\alpha}=f$ for all $\alpha\in\tilde{V}_{k}$ then $D_{k}(f_{\alpha}:\alpha\in\tilde{V}_{k})(x) = D_{k}f(x)$ is a continuous function. More generally, by switching $x$ to $y$ one by one, one has,
\begin{multline}\label{eq 5.1}
|D_{k}(f_{\alpha}:\alpha\in\tilde{V}_{k})(x)-D_{k}(f_{\alpha}:\alpha\in\tilde{V}_{k})(y)|\leq\int|f_{e_1}(x+h_1)-f_{e_1}(y+h_1)|\prod_{\alpha\in\tilde{V}_{k};\alpha\not=e_1}|f_{\alpha}|(x+\alpha\cdot\vec{h})\,d\vec{h}+\cdots\\+\int|f_{\vec{1}}(x+h_1+\cdots+h_{k})-f_{\vec{1}}(y+h_1+\cdots+h_{k})|\prod_{\alpha\in\tilde{V}_{k};\alpha\not=\vec{1}}|f_{\alpha}|(y+\alpha\cdot\vec{h})\,d\vec{h}.\end{multline}
Here $e_1=(1,0,\cdots,0),\vec{1}=(1,1,\cdots,1)\in\{0,1\}^{k}$. Define new functions $g_{\beta,\gamma}$ based on translations of $f_{\alpha}$ in order to utilize $\ref{eq 2.1}$ - here, $\alpha,\beta,\gamma\in\tilde{V}_{k}$. For instance, for the first factor of the first term on the RHS of $\ref{eq 5.1}$, one can define $g_{e_1,e_1}(z) = f_{e_1}(z)-f^{y-x}_{e_1}(z)$, then $g_{e_1,e_1}(x+h_1) = f_{e_1}(x+h_1) - f_{e_1}(y+h_1)$. With these new $g_{\beta,\gamma}$, by $\ref{eq 2.1}$, the LHS of $\ref{eq 5.1}$ is then dominated by,
\begin{equation}\label{eq 5.2}
|D_{k}(f_{\alpha}:\alpha\in\tilde{V}_{k})(x)-D_{k}(f_{\alpha}:\alpha\in\tilde{V}_{k})(y)|\lesssim\sum_{\beta\in\tilde{V}_{k}}\prod_{\gamma\in\tilde{V}_{k}}\|g_{\beta,\gamma}\|_{q_{k}}.\end{equation}
For each of the products on the RHS of $\ref{eq 5.2}$, only one factor in it takes the form $\|f_{\gamma}(\cdot)-f_{\gamma}^{y-x}(\cdot)\|_{q_{k}}$, for some $\gamma\in\tilde{V}_{k}$; each of the remaining factors is equal to $\|f_{\alpha}\|_{q_{k}}$. The conclusion then follows from the continuity of $L^{q_{k}}$ norms with respect to translation. $\qed$\\

\noindent
A consequence of {\bf Lemma 7} and $\ref{eq 2.1}$ is that, $\mathcal{K}(k)$ is contained in the closed unit ball of continuous functions with respect to the uniform norm - in fact, each function of $\mathcal{K}(k)$ belongs to the Banach algebra of bounded continuous functions with the uniform norm. Much like that the continuous functions form an algebra $\mathcal{C}$, one can also form an algebra out of $\mathcal{K}(k)$. Let $\bar{\mathcal{K}}(k)$ denote the closure of $\mathcal{K}(k)$ in $(\mathcal{C},\|\cdot\|_{\infty})$ and let $\mathbb{A}(k)$ be the linear subspace spanned by $\bar{\mathcal{K}}(k)$ endowed with the norm $\|\cdot\|_{\mathbb{A}(k)}$ such that $\bar{\mathcal{K}}(k)$ is the unit ball of $\mathbb{A}(k)$. More precisely, if $g\in\bar{\mathcal{K}}(k)$ is the limit of a non-constant sequence in $\mathcal{K}(k)$ under the uniform norm, then $\|g\|_{\mathbb{A}(k)}=1$. \\

\noindent
{\bf Lemma 8.} $(\mathbb{A}(k),\|\cdot\|_{\mathbb{A}(k)})$ is a Banach space.\\
{\it Proof:} Firstly, note that if $g\in\mathbb{A}(k)$ then \begin{equation}\label{eq 5.3}
\|g\|_{\infty}\leq\|g\|_{\mathbb{A}(k)}.\end{equation} 
Indeed, if $g\in\mathbb{A}(k)\cap\bar{\mathcal{K}}(k)$ such that $\|g\|_{\mathbb{A}(k)}= 1$, then from $\ref{eq 2.1}$ and the definition of $\bar{\mathcal{K}}(k)$, $\|g\|_{\infty}\leq 1=\|g\|_{\mathbb{A}(k)}$. The general statement then follows from the definition of $\mathbb{A}(k)$.\\

\noindent
Let $g_{n}$ be a sequence in $\mathbb{A}(k)$ such that $\sup_{n}\|g_{n}\|_{\mathbb{A}(k)}<\infty$ and that $g_{n}$ converges uniformly to some $g\in\mathcal{C}$. Then $g\in\mathbb{A}(k)$ and $\|g\|_{\mathbb{A}(k)}\leq\sup_{n}\|g_{n}\|_{\mathbb{A}(k)}<\infty$. This fact also follows from the definition of $\mathbb{A}(k)$. \\

\noindent
Now let $h_{n}\in\mathbb{A}(k)$ such that $\sum_{n}\|h_{n}\|_{\mathbb{A}(k)}<\infty$. From $\ref{eq 5.3}$, $\sum_{n}h_{n}$ converges uniformly to some $h=_{\|\cdot\|_{\infty}}\sum_{n}h_{n}$. Moreover, $h\in\mathbb{A}(k)$ with $\|h\|_{\mathbb{A}(k)}\leq\sum_{n}\|h_{n}\|_{\mathbb{A}(k)}$ by the observation above. It's easily seen that $\|h-\sum_{n=1}^{N}h_{n}\|_{\mathbb{A}(k)}\leq\sum_{n>N}\|h_{n}\|_{\mathbb{A}(k)}\to 0$ as $N\to\infty$.\\

\noindent
{\bf Theorem 9.} $\mathbb{A}(k)$ is a Banach algebra.\\
{\it Proof:} One must show that if $g^1, g^2\in\mathbb{A}(k)$ then so is $g^1g^2$ and $\|g^1g^2\|_{\mathbb{A}(k)}\leq\|g^1\|_{\mathbb{A}(k)}\|g^2\|_{\mathbb{A}(k)}$. By density, $\ref{eq 2.1}$ and the definition of $\mathbb{A}(k)$, it's enough to show that, if $g^{i}=D_{k}(f^{i}_{\alpha}:\alpha\in\tilde{V}_{k})$ with $f^{i}_{\alpha}\in\mathcal{S}$ and $\|f^{i}_{\alpha}\|_{q_{k}}\leq 1$, $i=1,2$, then $g^1g^2\in\mathbb{A}(k)$ and $\|g^1g^2\|_{\mathbb{A}(k)}\leq 1$.\\
Now, 
\begin{multline}\label{eq 5.4} g^1(x)g^2(x)=\int\left(\int\prod_{\alpha\in\tilde{V}_{k}}f^1_{\alpha}(x+\alpha\cdot\vec{h})f^2_{\alpha}(x+\alpha\cdot\vec{t})\,d\vec{h}\right)\,d\vec{t} \\= \int\left(\int\prod_{\alpha\in\tilde{V}_{k}}f^1_{\alpha}(x+\alpha\cdot\vec{h})f^2_{\alpha}(x+\alpha\cdot\vec{h}+\alpha\cdot\vec{u})\,d\vec{h}\right)\,d\vec{u}\end{multline}
with $\vec{u}=\vec{t}-\vec{h}$. Abstract way the dependence on $x$ for notational convenience. Integrating the last integral expression in $\ref{eq 5.4}$ first in terms of $\vec{u}$ and second in terms of $\vec{h}$ - and utilizing $\ref{eq 2.1}$ step by step - one has
\begin{equation}\label{eq 5.5} 
|g^1g^2|\lesssim\int \prod_{\alpha\in\tilde{V}_{k}}|f_{\alpha}^1|(\alpha\cdot\vec{h})\left(\prod_{\alpha\in\tilde{V}_{k}}\|f^2_{\alpha}\|_{q_{k}}\right)\,d\vec{h}\lesssim\prod_{\alpha\in\tilde{V}_{k}}\|f_{\alpha}^1\|_{q_{k}}\|f_{\alpha}^2\|_{q_{k}}\lesssim 1.\end{equation}
The constant involved in the last inequality of $\ref{eq 5.5}$ are multiples of the constants in sharp Young's inequality (see the proof of {\bf Lemma 1}) which are all less than $1$. Hence, 
\begin{equation}\label{eq 5.6}
|g^1(x)g^2(x)|\lesssim\prod_{\alpha\in\tilde{V}_{k}}\|f^1_{\alpha}\|_{q_{k}}\|f^2_{\alpha}\|_{q_{k}}\leq 1.\end{equation}
Since $f^1_{\alpha},f^2_{\alpha}\in\mathcal{S}$, a similar argument as in the proof of {\bf Lemma 7} gives that $g^1g^2$ is a continuous function. Then $\ref{eq 5.6}$ gives that $g^1g^2$ belongs to the Banach algebra of bounded continuous functions with the uniform norm. Now by definition, and by density, $g^1$ is approximated in the uniform norm by a sequence of the form $\{D_{k}(f^{1,m}_{\alpha}:\alpha\in\tilde{V}_{k})\}_{m}$ and $g^2$ by $\{D_{k}(f^{2,n}_{\alpha}:\alpha\in\tilde{V}_{k})\}_{n}$, with $f^{1,m}_{\alpha},f^{2,n}_{\alpha}\in\mathcal{S}$. Then $g^1g^2$ is approximated in the uniform norm by $\{D_{k}(f^{1,l}_{\alpha}:\alpha\in\tilde{V}_{k})\cdot D_{k}(f^{2,l}_{\alpha}:\alpha\in\tilde{V}_{k})\}_{l}$. Each term in the mentioned sequence has a form similarly to the last expression of $\ref{eq 5.4}$:
\begin{equation*}D_{k}(f^{1,l}_{\alpha}:\alpha\in\tilde{V}_{k})\cdot D_{k}(f^{2,l}_{\alpha}:\alpha\in\tilde{V}_{k})(x) = \int\left(\int\prod_{\alpha\in\tilde{V}_{k}}f^{1,l}_{\alpha}(x+\alpha\cdot\vec{h})f^{2,l}_{\alpha}(x+\alpha\cdot\vec{h}+\alpha\cdot\vec{u})\,d\vec{h}\right)\,d\vec{u}.\end{equation*}
By the Schwartz nature, one can approximate this last integral expression by a Stieljes sum, 
\begin{multline*}C\sum_{i=1}^{N}\lambda_{i}\int\prod_{\alpha\in\tilde{V}_{k}}f^{1,l}_{\alpha}(x+\alpha\cdot\vec{h})f^{2,l}_{\alpha}(x+\alpha\cdot\vec{h}+\alpha\cdot\vec{u}_{i})\,d\vec{h}\\=\sum_{i=1}^{N}\lambda_{i}\int C^{1/(2^{k}-1)}\prod_{\alpha\in\tilde{V}_{k}}f^{1,l}_{\alpha}(x+\alpha\cdot\vec{h})f^{2,l}_{\alpha}(x+\alpha\cdot\vec{h}+\alpha\cdot\vec{u}_{i})\,d\vec{h} =: \sum_{i=1}^{N}\lambda_{i} D_{k}(w^{l}_{\alpha}:\alpha\in\tilde{V}_{k})\end{multline*}
where $\lambda_{i}>0$ and $\sum_{i}\lambda_{i}=1$ and $C$ is a normalized constant ($C\lambda_{i}=|\Delta_{i}|$, the discretization area units). Again, by the Schwartz nature, $w^{l}_{\alpha}\in L^{q_{k}}\cap L^{p_{k}}$. Summarizing all of this, one has that for some small $\delta$,
\begin{equation*}
\sup_{x}|g^1(x)g^2(x)-\sum_{i=1}^{N}\lambda_{i} D_{k}(w^{l}_{\alpha}:\alpha\in\tilde{V}_{k})(x)|<\delta
\end{equation*}
Then from $\ref{eq 5.6}$, 
\begin{equation*}\sup_{x}|\sum_{i=1}^{N}\lambda_{i} D_{k}(w^{l}_{\alpha}:\alpha\in\tilde{V}_{k})(x)|<1+\delta,\end{equation*}
which means $\sum_{i=1}^{N}\lambda_{i} D_{k}(w^{l}_{\alpha}:\alpha\in\tilde{V}_{k})\in (1+\delta)\bar{\mathcal{K}}$, following from definition. That means $g^1g^2$ is approximated in the uniform norm by a function in $(1+\delta)\bar{\mathcal{K}}$ for every $\delta$. Let $\delta\to 0$, one concludes that $g^1g^2\in\bar{\mathcal{K}}$. Then $\|g^1g^2\|_{\mathbb{A}(k)}\leq 1$. $\qed$

\subsection{The algebra $\mathbb{F}(k)$}

\noindent
In the light of {\bf Theorem 4}, suppose instead of $\mathcal{K}(k)$, one considers $\mathcal{L}(k)$ that is the convex hull of
\begin{equation*}\{D_{k}(f_{\alpha}:\alpha\in\tilde{V}_{k}): f_{\alpha}\in L^{q_{k}}\cap L^{p_{k}}\text{ and } \prod_{\alpha\in\tilde{V}_{k}}\|f_{\alpha}\|_{q_{k}}\leq 1\}.\end{equation*}
Then {\bf Lemma 7} and $\ref{eq 2.1}$ still hold. Then similarly, let $\bar{\mathcal{L}}(k)$ denote the closure of $\mathcal{L}(k)$ in $(\mathcal{C},\|\cdot\|_{\infty})$ and let $\mathbb{F}(k)$ be the linear subspace spanned by $\bar{\mathcal{L}}(k)$ endowed with the norm $\|\cdot\|_{\mathbb{F}(k)}$ such that $\bar{\mathcal{K}}(k)$ is the unit ball of $\mathbb{F}(k)$. An analogous version of {\bf Lemma 8} and its proof holds for the case of $\mathbb{F}(k)$:\\

\noindent
{\bf Lemma 10.} $(\mathbb{F}(k),\|\cdot\|_{\mathbb{F}(k)})$ is a Banach space.\\

\noindent
Since the proof of {\bf Theorem 9} uses approximation by Schwartz functions, it's still carried over to the context of $\mathbb{F}(k)$.\\

\noindent
{\bf Theorem 11.} $\mathbb{F}(k)$ is a Banach algebra.\\

\noindent
As the matter of fact, the reason to consider functions in $L^{p_{k}}\cap L^{q_{k}}$ in the definition of $\mathbb{F}(k)$ is the following.\\

\noindent
Let $g(x) = D_{k}(f_{\alpha}:\alpha\in\tilde{V}_{k})(x)$ with $f_{\alpha}\in L^{p_{k}}\cap L^{q_{k}}$. Follow the same argument as in $\ref{eq 3.1}$ (or use $\ref{eq 1.5}$) and utilize $\ref{eq 2.1}$, one can deduce $g\in L^{r}, s_{k}\leq r\leq\infty$, and in particular, such $g$ belongs to the Banach algebra of continuous, bounded functions that vanish at infinity. Moreover,
\begin{equation*}\|g\|_{s_{k}}\lesssim\prod_{\alpha\in\tilde{V}_{k}}\|f_{\alpha}\|_{p_{k}}.\end{equation*}
It's easily seen that $p_{k}$'s are increasing as $k\uparrow\infty$. If $k\geq 3$ then $p_{k}\geq 2$, hence $1\leq s_{k}\leq 2$. From the boundedness of Fourier transform on $L^{r}$, $1\leq r\leq 2$, one has
\begin{equation}\label{eq 5.7}\|\hat{g}\|_{p_{k}}\lesssim\prod_{\alpha\in\tilde{V}_{k}}\|f_{\alpha}\|_{p_{k}}.\end{equation}
When $k=3$, $p_{k}=2$ and $\ref{eq 5.7}$ becomes, $\|\hat{g}\|_2\lesssim\prod_{\alpha\in\tilde{V}_{k}}\|f_{\alpha}\|_2$. Using density, $\ref{eq 5.7}$ can be extended to functions $g=D_{k}(f_{\alpha}:\alpha\in\tilde{V}_{k})$ with $f_{\alpha}\in L^{q_{k}}$, $k\geq 3$:
\begin{equation*}\|\hat{g}\|_{p_{k}}\lesssim\prod_{\alpha\in\tilde{V}_{k}}\|f_{\alpha}\|_{q_{k}}.\end{equation*}

\section{Remarks}

\noindent
As the Brascamp-Lieb-type inequalities have a strong connection to geometric functional analysis on Euclidean spaces, see [2], [3], [4], [10], [9], [11], the author hopes that this initial exploration sets up a longer investigation into the world of the anti-uniform spaces in the future and their geometric implications.



\begin{thebibliography}{}

\bibitem{Br}
H. Brezis,
\textit{Functional Analysis, Sobolev spaces, and partial differential equations},
Universitext, Springer, 2010.

\bibitem{Bu}
A. Burchard,
\textit{Cases of equality in the Riesz rearrangement inequality},
Ann. of Math., (2) 143 (1996), no. 3, 499-528.

\bibitem{BCChT}
J. Bennett, A. Carbery, M. Christ, T. Tao,
\textit{The Brascamp-Lieb inequalities: finiteness, structure, and extremals},
Geom. Funct. Anal., 17 (2008), no. 5, 1343-1415.

\bibitem{Ch}
M. Christ,
\textit{Subsets of Euclidean space with nearly maximal Gowers norms},
preprint, math.CA arXiv:1512.03355.

\bibitem{ET}
T. Eisner and T. Tao,
\textit{Large values of the Gowers-Host-Kra seminorms},
J. Anal. Math., 117 (2012), 133-186. 

\bibitem{G}
T. Gowers,
\textit{A new proof of Szemer\'edi's theorem},
Geom. Funct. Anal., 11 (2001), 465-588.

\bibitem{HK1}
B. Host and B. Kra,
\textit{Nonconventional averages and nilmanifolds},
Ann. of Math., 161 (2005) 398-488.

\bibitem{HK2}
B. Host and B. Kra,
\textit{A point of view on Gowers uniformity norms},
New York J. Math., 18 (2010).

\bibitem{L}
E. Lieb, 
\textit{Gaussian kernels have only Gaussian maximizers},
J.L. Invent. Math., 102 (1990), 102-179.

\bibitem{LL}
E. Lieb and M. Loss, 
\textit{Analysis},
AMS, Rrovidence, RI, (2001).

\bibitem{N}
M. Neuman,
\textit{Functions of Nearly Maximal Gowers–Host–Kra Norms on Euclidean Spaces},
J. Geom. Anal. (2018), 1-58.

\bibitem{combinatorics}
T. Tao and V. Vu,
\textit{Additive Combinatorics},
Cambridge University Press, 2016.

\end{thebibliography}
\end{document}